\documentclass{article}
\usepackage{amsmath}
\usepackage{amsfonts}
\usepackage{amssymb}
\usepackage{graphicx}
\usepackage{cite}

\newtheorem{theorem}{Theorem}[section]

\newtheorem{definition}[theorem]{Definition}
\newtheorem{example}[theorem]{Example}

\newtheorem{lemma}[theorem]{Lemma}

\input{tcilatex}
\begin{document}

\title{On the Numerical Solution of Nonlinear Fractional-Integro
Differential Equations}
\author{Mehmet \c{S}ENOL and \.{I}. Timu\c{c}in DOLAPC\.{I} \\
Nev\c{s}ehir Hac\i\ Bekta\c{s} Veli University, Department of Mathematics,
Nev\c{s}ehir, Turkey\\
Celal Bayar University, Department of Mechanical Engineering,\\
Manisa, Turkey\\
e-mail:msenol@nevsehir.edu.tr, ihsan.dolapci@cbu.edu.tr}
\maketitle

\begin{abstract}
In the present study, a numerical method, perturbation-iteration algorithm
(shortly PIA), have been employed to give approximate solutions of nonlinear
fractional-integro differential equations (FIDEs). Comparing with the exact
solution, the PIA produces reliable and accurate results for FIDEs.

\end{abstract}

\section{Introduction}

Scientists has been interested in fractional order calculus as long as it has been in classical integer order analysis.
However, for many years it could not find practical applications in physical sciences. Recently, fractional calculus has been used in applied
mathematics, viscoelasticity \cite{1}, control \cite{2}, electrochemistry 
\cite{3}, electromagnetic \cite{4}. 

Developments in symbolic computation capabilities is one of the driving forces behind this rise. Different multidisciplinary problems can be handled with fractional derivatives and integrals.

\cite{5} and \cite{6} are studies that describe the fundamentals of fractional calculus give some applications. Existence and uniqueness of the solutions
are also studied in \cite{7}.

Similar to the studies in physical sciences, fractional order integro differential equations (FIDEs) also gave scientists the opportunity of describing and modeling many important and useful physical problems.

In this manner, a remarkable effort has been expended to propose numerical methods for solving FIDEs, in recent years. Fractional
variational iteration method \cite{8,9}, homotopy analysis method \cite%
{10,11}, Adomian decomposition method \cite{12,13} and fractional
differential transform method \cite{14,15,16} are among these methods.

In our study, we use the previously developed method PIA, to obtain approximate solutions of some FIDEs. This method can be applied to a wide range of problems without requiring any special assumptions and restrictions.

A few fractional derivative definitions of an arbitrary order exists in the literature.  Two most used of them are the Riemann-Liouville and Caputo fractional derivatives. The two definitions are quite similar but have different order of evaluation of derivation. 

\qquad The Riemann-Liouville fractional integral of order $\alpha $ is
described by:%
\begin{equation}
J^{\alpha }u(x)=\frac{1}{\Gamma (\alpha )}\int_{0}^{x}(x-t)^{\alpha
-1}u(t)dt,\quad \alpha >0,\quad x>0.  \label{1}
\end{equation}

\qquad The Riemann-Liouville and Caputo fractional derivatives of an
arbitrary order are defined as the following, respectively%
\begin{equation}
D^{\alpha }u(x)=\frac{d^{m}}{dx^{m}}\left( J^{m-\alpha }u(x)\right)
\label{2}
\end{equation}%
\begin{equation}
D_{\ast }^{\alpha }u(x)=J^{m-\alpha }\left( \frac{d^{m}}{dx^{m}}u(x)\right) .
\label{3}
\end{equation}%
where $m-1<\alpha \leqslant m$ and $m\in 
\mathbb{N}
.$

Due to the appropriateness of the initial conditions, fractional definition
of Caputo is often used in recent years.

\begin{definition}
The Caputo fractional derivative of a function $u(x)$ is defined as%
\begin{equation}
D_{\ast }^{\alpha }u(x)=\left\{ 
\begin{array}{cc}
\frac{1}{\Gamma (m-\alpha )}\int_{0}^{x}(x-t)^{m-\alpha -1}u^{(m)}(t)dt, & 
m-1<\alpha \leqslant m \\ 
\frac{d^{m}u(x)}{dx^{m}} & \alpha =m%
\end{array}%
\right.  \label{4}
\end{equation}%
for $m-1<\alpha \leqslant m,$ $m\in 
\mathbb{N}
,$ $x>0,$ $u\in C_{-1}^{m}.$
\end{definition}

Following lemma gives the two main properties of Caputo fractional
derivative.

\begin{lemma}
For $m-1<\alpha \leqslant m,$ $u\in C_{\mu }^{m},$ $\mu \geqslant -1$ and $%
m\in 
\mathbb{N}
$ then 
\begin{equation}
D_{\ast }^{\alpha }J^{\alpha }u(x)=u(x)  \label{5}
\end{equation}%
and 
\begin{equation}
J^{\alpha }D_{\ast }^{\alpha }u(x)=u(x)-\sum_{k=0}^{m-1}u^{(k)}(0^{+})\frac{%
x^{k}}{k!},\quad x>0.  \label{6}
\end{equation}
\end{lemma}

After this introductory section, Section 2 is reserved to a brief review of
the Perturbation-Iteration Algorithm PIA, in Section 3 some examples are
illustrated to show the simplicity and effectiveness of the algorithm.
Finally the paper ends with a conclusion in Section 4.

\section{Analysis of the PIA}

Differential equations are naturally used to describe problems in engineering and other applied sciences. Searching approximate solutions for complicated equations has always attracted attention. Many different methods and frameworks exist for this purpose and perturbation techniques \cite{17,18,19} are among them. These techniques can be used to find approximate solutions for both ordinary and partial differential equations. 

Requirement of a small parameter in the equation that is sometimes artificially inserted is a major limitation of perturbation techniques that renders them valid only in a limited range. Therefore, to overcome the disadvantages come with the perturbation techniques, several methods have been proposed by authors \cite{20,21,22,23,24,25,26,27,28,29}.

Parallel to these attempts, a perturbation-iteration method has been
proposed by Aksoy, Pak\-de\-mir\-li and their co-workers \cite{33,34,35}
previously. A primary effort of producing root finding algorithms for
algebraic equations \cite{30,31,32}, finally guided to obtain formulae for
differential equations also \cite{33,34,35}. In the new technique, an
iterative algorithm is constructed on the perturbation expansion. The
present method has been tested on Bratu-type differential equations \cite{33}
and first order differential equations \cite{34} with success. Then the
algorithms were applied to nonlinear heat equations also \cite{35}. Finally,
the solutions of the Volterra and Fredholm type integral equations \cite{36}
and ordinary differential equation systems \cite{37} have been presented by
the developed method.

This new algorithm have not been used for any fractional integro differential equations
yet. To obtain the approximate solutions of FIDEs, the most basic
perturbation-iteration algorithm PIA(1,1) is employed by taking one
correction term in the perturbation expansion and correction terms of only
first derivatives in the Taylor series expansion. \cite{33,34,35}.

Take the fractional-integro differential equation. 
\begin{equation}
F\left( u^{(\alpha )},u,\int_{0}^{t}{g\left( t,s,u(s)\right) ds},\varepsilon
\right) =0  \label{7}
\end{equation}

where $u=u(t)$ and $\varepsilon $ is a small parameter. The perturbation
expansions with only one correction term is

\begin{eqnarray}
u_{n+1} &=&u_{n}+\varepsilon {\left( u_{c}\right) }_{n}\   \notag \\
u_{n+1}^{\prime } &=&u_{n}^{\prime }+\varepsilon {\left( u_{c}^{\prime
}\right) }_{n}\   \label{8}
\end{eqnarray}

Replacing Eq.$(\ref{8})$ into Eq.$(\ref{7})$ and writing in the Taylor series
expansion for only first order derivatives gives

\begin{eqnarray}
&&F\left( u_{n}^{\left( \alpha \right) },u_{n},\int_{0}^{t}{g\left(
t,s,u_{n}(s)\right) ds},0\right)  \notag \\
&&+F_{u}\left( u_{n}^{\left( \alpha \right) },u_{n},\int_{0}^{t}{g\left(
t,s,u_{n}(s)\right) ds},0\right) \varepsilon {\left( u_{c}\right) }_{n} 
\notag \\
&&+F_{u^{\left( \alpha \right) }}\left( u_{n}^{\left( \alpha \right)
},u_{n},\int_{0}^{t}{g\left( t,s,u_{n}(s)\right) ds},0\right) \varepsilon {%
\left( u_{c}^{(\alpha )}\right) }_{n}  \notag \\
&&+F_{\int {u}}\left( u_{n}^{\left( \alpha \right) },u_{n},\int_{0}^{t}{%
g\left( t,s,u_{n}(s)\right) ds},0\right) \varepsilon \int {{\left(
u_{c}\right) }_{n}}  \notag \\
&&+F_{\varepsilon }\left( u_{n}^{\left( \alpha \right) },u_{n},\int_{0}^{t}{%
g\left( t,s,u_{n}(s)\right) ds},0\right) \varepsilon =0  \label{9}
\end{eqnarray}

or

\begin{equation}
{\left( u_{c}^{(\alpha )}\right) }_{n}\frac{\partial F}{\partial u^{(\alpha
)}}+{\left( u_{c}\right) }_{n}\frac{\partial F}{\partial u}+\left( \int {{%
\left( u_{c}\right) }_{n}}\right) \frac{\partial F}{\partial (\int {u})}+%
\frac{\partial F}{\partial \varepsilon }+\frac{F}{\varepsilon }=0  \label{10}
\end{equation}

Here $(.)^{\prime }$ represents the derivative according to the independent
variable and 
\begin{equation}
F_{\varepsilon }=\frac{\partial F}{\partial \varepsilon },~F_{u}=\frac{%
\partial F}{\partial u},~F_{u^{\prime }}=\frac{\partial F}{\partial
u^{\prime }},\ldots  \label{11}
\end{equation}%

The derivatives in the expansion are evaluated at $\varepsilon =0$.
Beginning with an initial function $u_{0}(t)$, first ${\left( u_{c}\right) }%
_{0}(t)$ is calculated by the help of $(\ref{10})$ and then substituted into Eq.$(\ref{8})$ to calculate $%
u_{1}(t)$. Iteration procedure is continued using $(\ref{10}) $
and $(\ref{8})$ until obtaining a
reasonable solution.

\section{Applications}

\begin{example}
Consider the following nonlinear fractional-integro differential equation 
\cite{38}:
\end{example}

\begin{equation}
\frac{d^{\alpha }u(t)}{{dt}^{\alpha }}-\int_{0}^{1}{ts({u(s))}^{2}ds}=1-%
\frac{t}{4},\ \ \ t>0,\ \ \ 0\leq t<1,\ \ \ 0<\alpha \leq 1  \label{12}
\end{equation}

with the initial condition\ $u\left( 0\right) =0$ and the known exact
solution for $\alpha =1$ is

\begin{equation}
u\left( t\right) =t  \label{13}
\end{equation}

Before iteration process rewriting Eq.$(\ref{12})$ with adding and
subtracting $u^{\prime }(t)$ to the equation gives

\begin{equation}
\varepsilon \frac{d^{\alpha }u(t)}{{dt}^{\alpha }}-u^{^{{\prime }}}\left(
t\right) +{\varepsilon u}^{^{{\prime }}}\left( t\right) -{\varepsilon }%
\int_{0}^{1}{ts({u(s))}^{2}ds}-1+\frac{t}{4}=0  \label{14}
\end{equation}%
In this case for

\begin{equation}
F\left( u^{^{{\prime }}},u,\varepsilon \right) =\frac{1}{\Gamma (1-\alpha )}%
\varepsilon \int_{0}^{t}{\frac{u^{\prime }(s)}{{(t-s)}^{\alpha }}ds-u_{n}^{^{%
{\prime }}}\left( t\right) +\varepsilon u_{n}^{^{{\prime }}}\left( t\right)
-\varepsilon \int_{0}^{1}{ts({u_{n}(s))}^{2}ds}-1+\frac{t}{4}}  \label{15}
\end{equation}%
and the iteration formula%
\begin{equation}
u^{^{{\prime }}}(t)+\frac{F_{u}}{F_{u^{\prime }}}u\left( t\right) =-\frac{%
F_{\varepsilon }+\frac{F}{\varepsilon }}{F_{u^{\prime }}}  \label{16}
\end{equation}

the terms that will be replaced in, are

\begin{eqnarray}
F &=&{u_{n}^{^{{\prime }}}\left( t\right) }-1+\frac{t}{4}  \notag \\
F_{u} &=&0  \notag \\
F_{u^{\prime }} &=&1  \notag \\
F_{\varepsilon } &=&-{u_{n}^{^{{\prime }}}\left( t\right) }+\frac{1}{\Gamma
(1-\alpha )}\int_{0}^{t}{\frac{u^{\prime }(s)}{{(t-s)}^{\alpha }}ds}%
-\int_{0}^{1}{ts({u(s))}^{2}ds}  \label{17}
\end{eqnarray}

After substitution the differential equation for this problem, Eq.$(\ref%
{10})$ becomes \bigskip

\begin{equation}
\frac{\int_{0}^{t}{{\left( -s+t\right) }^{-\alpha }{u_{n}}^{^{{\prime }%
}}(s)ds}}{\Gamma (1-\alpha )}+{\left( {u}_{c}^{\prime }(t)\right) }%
_{n}=\int_{0}^{1}{st{\left( u_{n}\left( s\right) \right) }^{2}ds}+\frac{%
4-t+4\left( -1+\varepsilon \right) {u}_{n}^{\prime }\left( t\right) }{%
4\varepsilon }  \label{18}
\end{equation}

Appropriate to the initial conditions, chosen $u_{0}\left( t\right) =0$ and,
solving Eq.$(\ref{18})$ for $n=0$ gives 
\begin{equation}
{{(u}_{c}(t))}_{0}=t-\frac{t^{2}}{8}+C_{1}  \label{19}
\end{equation}

This expression written in

\begin{equation}
u_{1}=u_{0}+\varepsilon {{(u}_{c}(t))}_{0}  \label{20}
\end{equation}

gives

\begin{equation}
u_{1}\left( x,t\right) =u_{0}\left( x,t\right) +\varepsilon (t-\frac{t^{2}}{8%
}+C_{1})  \label{21}
\end{equation}

or

\begin{equation}
u_{1}\left( x,t\right) =\varepsilon (t-\frac{t^{2}}{8}+C_{1})  \label{22}
\end{equation}

Solving this equation for

\begin{equation}
u_{1}\left( 0\right) =0  \label{23}
\end{equation}

we obtain

\begin{equation}
C_{1}=0  \label{24}
\end{equation}%
For this value and $\varepsilon =1$ reorganizing $u_{1}(t)$

\begin{equation}
u_{1}\left( t\right) =t-\frac{t^{2}}{8}  \label{25}
\end{equation}

gives the first iteration result. If the iteration procedure is continued in
a similar way, we obtain the following iterations.

\begin{equation}
u_{2}(t)=2t-\frac{571t^{2}}{3840}+\frac{t^{2-\alpha }(t+4(-3+\alpha ))}{%
4\Gamma (4-\alpha )}\   \label{26}
\end{equation}%
\begin{eqnarray}
u_{3}\left( t\right) &=&3t+\frac{29844889t^{2}}{176947200}-\frac{%
t^{3-2\alpha }\left( t+8\left( -2+\alpha \right) \right) }{4\Gamma \left(
5-2\alpha \right) }  \notag \\
&&+\frac{t^{2}\left( 3379230+8t^{-\alpha }\left( 1051t+5760\left( -3+\alpha
\right) \right) \left( -7+\alpha \right) \left( -6+\alpha \right) \left(
-5+\alpha \right) \right) }{15360\left( -7+\alpha \right) \left( -6+\alpha
\right) \left( -5+\alpha \right) \Gamma \left( 4-\alpha \right) }  \notag \\
&&-\frac{2240277\alpha +\left( 450151-28436\alpha \right) \alpha ^{2}}{%
15360\left( -7+\alpha \right) \left( -6+\alpha \right) \left( -5+\alpha
\right) \Gamma \left( 4-\alpha \right) }  \notag \\
&&-\frac{t^{2}\left( -4+\alpha \right) \left( -1159+2\alpha \left(
529+16\left( -10+\alpha \right) \alpha \right) \right) }{64\left( -7+2\alpha
\right) {\Gamma \left( 5-\alpha \right) }^{2}}  \label{27}
\end{eqnarray}

The other iterations contain large inputs and are not given. A computational
software program could help to calculate the other iterations up to any
order. In Table 1. some of the PIA iteration results are compared with
the exact solution. The results express that the present method gives highly approximate solutions. Also in Figure 1. the obtained results are
illustrated graphically.

\begin{table}[th]
\caption{Numerical results of Example 3.1. for different $u$ values when $%
\protect\alpha =1$ }
\begin{center}
{\scriptsize 
\begin{tabular}{ccccccc}
\hline
& \multicolumn{6}{c}{$\alpha =1.0$} \\ \hline
t & $u_{2}$ & $u_{3}$ & $u_{4}$ & $u_{5}$ & Exact Solution & Absolute Error
\\ \hline
$0.0$ & 0.000000 & 0.000000 & 0.000000 & 0.000000 & 0.000000 & 0.000000 \\ 
$0.1$ & 0.099763 & 0.099953 & 0.099990 & 0.099981 & 0.100000 & 1.872712E-6
\\ 
$0.2$ & 0.199052 & 0.199812 & 0.199962 & 0.199992 & 0.200000 & 7.490848E-6
\\ 
$0.3$ & 0.297867 & 0.299577 & 0.299915 & 0.299983 & 0.300000 & 1.685440E-5
\\ 
$0.4$ & 0.396208 & 0.399249 & 0.399850 & 0.399970 & 0.400000 & 2.996339E-5
\\ 
$0.5$ & 0.494075 & 0.498826 & 0.499765 & 0.499953 & 0.500000 & 4.681780E-5
\\ 
$0.6$ & 0.591468 & 0.598310 & 0.599662 & 0.599932 & 0.600000 & 6.741763E-5
\\ 
$0.7$ & 0.688388 & 0.697700 & 0.699541 & 0.699908 & 0.700000 & 9.176289E-5
\\ 
$0.8$ & 0.784833 & 0.796996 & 0.799400 & 0.799880 & 0.800000 & 1.198535E-4
\\ 
$0.9$ & 0.880804 & 0.896198 & 0.899241 & 0.899848 & 0.900000 & 1.516896E-4
\\ 
$1.0$ & 0.976302 & 0.995307 & 0.999063 & 0.999812 & 1.000000 & 1.872712E-4
\\ \hline
\end{tabular}
}
\end{center}
\end{table}

\begin{figure}[tbp]
\centering
\includegraphics[width=3.50in]{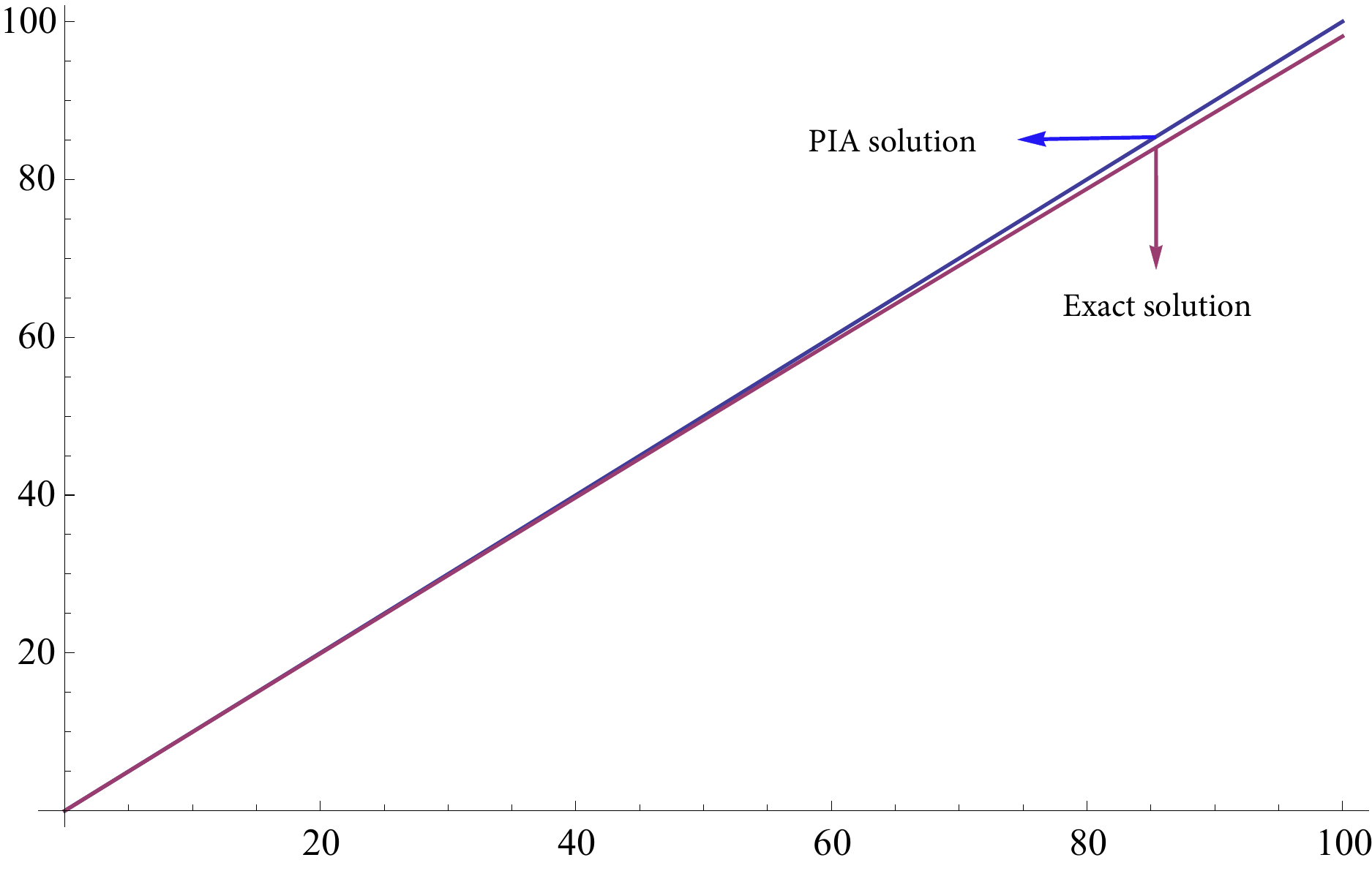}
\caption{Comparison of the PIA solution $u_{3}(t)$ and exact solution for
Example 3.1. when $\protect\alpha =1$}
\end{figure}

\begin{example}
Consider the following system of nonlinear fractional-integro differential
equations \cite{39}:
\end{example}

\begin{eqnarray}
\frac{d^{\alpha _{1}}u(t)}{{dt}^{\alpha _{1}}} &=&1-\frac{1}{2}{\left( k^{^{{%
\prime }}}\left( t\right) \right) }^{2}+\int_{0}^{t}{\left( \left(
t-s\right) k\left( s\right) +u\left( s\right) k\left( s\right) \right) ds} 
\notag \\
\frac{d^{\alpha _{2}}k(t)}{{dt}^{2}} &=&2t+\int_{0}^{t}{\left( \left(
t-s\right) u\left( s\right) -k^{2}\left( s\right) +u^{2}\left( s\right)
\right) ds}\ \ \ \ 0<\alpha _{1},\alpha _{2}\leq 1  \label{28}
\end{eqnarray}

Given with $u\left( 0\right) =0,\ \ k(0)=1$ as initial conditions. The exact
solution for $\alpha _{1} =\alpha _{2} =1$ is

\begin{eqnarray}
u\left( t\right) &=&sinht  \notag \\
k(t) &=&cosht  \label{29}
\end{eqnarray}%
Rewriting Eq.$(\ref{28})$ in the following for with adding and
subtracting $u^{\prime }(t)$ and $k^{\prime }(t)$ to the equation
respectively gives 
\begin{eqnarray}
&&\varepsilon \frac{d^{\alpha _{1}}u(t)}{{dt}^{\alpha _{1}}}+u^{\prime
}\left( t\right) -\varepsilon u^{\prime }(t)-1+\frac{1}{2}{\left( k^{\prime
}\left( t\right) \right) }^{2}-\varepsilon \int_{0}^{t}{\left( \left(
t-s\right) k\left( s\right) -u\left( s\right) k\left( s\right) \right) ds\ }
\notag \\
&&\varepsilon \frac{d^{\alpha _{2}}u(t)}{{dt}^{\alpha _{2}}}+k^{\prime
}\left( t\right) -\varepsilon k^{\prime }\left( t\right) -2t-\varepsilon
\int_{0}^{t}{\left( \left( t-s\right) u\left( s\right) +{k}^{2}\left(
s\right) -{u}^{2}\left( s\right) \right) ds}  \label{30}
\end{eqnarray}%
In this case for 
\begin{eqnarray}
F\left( u^{\prime },u,\varepsilon \right) &=&\frac{1}{\Gamma (1-\alpha _{1})}%
\varepsilon \int_{0}^{t}{\frac{u^{\prime }(s)}{{(t-s)}^{\alpha _{1}}}%
ds-\varepsilon \int_{0}^{t}{\left( \left( t-s\right) k\left( s\right)
+u\left( s\right) k\left( s\right) \right) ds}-1+\frac{1}{2}{\left(
k^{\prime }\left( t\right) \right) }^{2}}  \notag \\
F\left( k^{\prime },k,\varepsilon \right) &=&\frac{1}{\Gamma (1-\alpha _{2})}%
\varepsilon \int_{0}^{t}{\frac{u^{\prime }(s)}{{(t-s)}^{\alpha _{2}}}%
ds-\varepsilon \int_{0}^{t}{\left( \left( t-s\right) u\left( s\right)
-k^{2}\left( s\right) +u^{2}\left( s\right) \right) ds}}-2t  \label{31}
\end{eqnarray}

and the iteration formula

\begin{equation}
u^{\prime }\left( t\right) +\frac{Fu}{Fu^{\prime }}u\left( t\right) =-\frac{%
F_{\varepsilon }+\frac{F}{\varepsilon }}{Fu^{\prime }}  \label{32}
\end{equation}

the terms that will be replaced in, are

\begin{eqnarray}
F &=&u_{n}^{\prime }(t)-1+\frac{{k_{n}^{\prime }(t)}^{2}}{2}  \notag \\
F_{u} &=&0  \notag \\
F_{u^{\prime }} &=&1  \notag \\
F_{\varepsilon } &=&-u_{n}^{\prime }(t)+\frac{1}{\Gamma (1-\alpha _{1})}%
\int_{0}^{t}{\frac{u_{n}^{\prime }(s)}{{(t-s)}^{\alpha _{1}}}ds}-\int_{0}^{t}%
{((t-s)k_{n}(s)+u_{n}(s)k_{n}(s))ds}  \label{33}
\end{eqnarray}

\bigskip and the iteration formula%
\begin{equation}
k^{\prime }\left( t\right) +\frac{F_{k}}{F_{k^{\prime }}}k\left( t\right) =-%
\frac{F_{\varepsilon }+\frac{F}{\varepsilon }}{F_{k^{\prime }}}  \label{34}
\end{equation}

the terms that will be replaced in, are

\begin{eqnarray}
F &=&k_{n}^{\prime }\left( t\right) -2t\   \notag \\
F_{k} &=&0  \notag \\
F_{k^{\prime }} &=&1  \notag \\
F_{\varepsilon } &=&-k_{n}^{\prime }(t)+\frac{1}{\Gamma (1-\alpha _{2})}%
\int_{0}^{t}{\frac{k_{n}^{\prime }(s)}{{(t-s)}^{\alpha _{2}}}ds}-\int_{0}^{t}%
{((t-s)u_{n}(s)-{k_{n}(s)}^{2}+{u_{n}(s)}^{2})ds}  \label{35}
\end{eqnarray}

After substitution, the system of differential equations for this problem
become

\begin{equation*}
\frac{1}{\Gamma (1-\alpha _{1})}\int_{0}^{t}{{(-s+t)}^{-\alpha _{1}}{u}}%
^{\prime }{{_{n}}(s)ds}+{\left( {u^{\prime }}_{c}(t)\right) }_{n}+\frac{-1+%
\frac{1}{2}{k}^{\prime }{{_{n}}(t)}^{2}+{u}^{\prime }{_{n}}(t)}{\varepsilon }%
=\int_{0}^{t}{k_{n}(s)(-s+t+u_{n}(s))ds}+{u}^{\prime }{_{n}}(t)
\end{equation*}%
\begin{equation}
\frac{1}{\Gamma (1-\alpha _{2})}\int_{0}^{t}{{(-s+t)}^{-\alpha _{2}}{k}}%
^{\prime }{{_{n}}(s)ds}+{\left( {k^{\prime }}_{c}(t)\right) }%
_{n}=\int_{0}^{t}{(-{k_{n}(s)}^{2}+u_{n}(s)(-s+t+u_{n}(s)))ds}+\frac{%
2t+(-1+\varepsilon ){k}_{n}^{\prime }(t)}{\varepsilon }  \label{36}
\end{equation}

Appropriate to the initial conditions, chosen $u_{0}\left( t\right) =0$ and $%
k_{0}\left( t\right) =1$ and solving Eq.$(\ref{36})$ for $n=0,1,2,3,...$
the successive iterations are

\begin{equation}
u_{1}(t)=\frac{1}{6}(6t+t^{3})  \label{37}
\end{equation}

\begin{equation}
k_{1}(t)=1+\frac{t^{2}}{2}  \label{38}
\end{equation}

\begin{equation}
u_{2}\left( t\right) =\frac{1}{504}t\left(
1008+168t^{2}+21t^{4}+t^{6}\right) -\frac{t^{2-\alpha _{1}}\left(
12+t^{2}+\left( -7+\alpha _{1}\right) \alpha _{1}\right) }{\Gamma (5-\alpha
_{1})}  \label{39}
\end{equation}

\begin{equation}
k_{2}\left( t\right) =1+t^{2}+\frac{t^{4}}{24}+\frac{t^{6}}{240}+\frac{t^{8}%
}{2016}-\frac{t^{3-\alpha _{2}}}{\Gamma (4-\alpha _{2})}  \label{40}
\end{equation}

Following in this manner the third iteration results, $u_{3}(t)$ and $%
k_{3}(t),$ are calculated. Again Table 2, Figure 2 and Figure 3 prove
that PIA give remarkably approximate results. We can say that the higher iterations would
give closer results. \bigskip

\begin{table}[th]
\caption{Numerical results of Example 3.2. for $u_{3}$ and $k_{3}$ values
when $\alpha _{1} =\alpha _{2} =1$ }
\begin{center}
{\scriptsize 
\begin{tabular}{ccccccc}
\hline
& \multicolumn{6}{c}{$\alpha _{1} =\alpha _{2} =1$} \\ \hline
t & PIA $(u_{3})$ & Exact Solution & Absolute Error & PIA $(k_{3})$ & Exact
Solution & Absolute Error \\ \hline
$0.0$ & 0.000000 & 0.000000 & 0.000000 & 1.000000 & 1000000. & 0.000000 \\ 
$0.1$ & 0.100166 & 0.100166 & 1.591577E-10 & 1.005004 & 1.005004 & 
1.191735E-11 \\ 
$0.2$ & 0.201335 & 0.201336 & 2.053723E-8 & 1.020066 & 1.020066 & 3.060393E-9
\\ 
$0.3$ & 0.304519 & 0.304520 & 3.556439E-7 & 1.045338 & 1.045338 & 7.884730E-8
\\ 
$0.4$ & 0.410749 & 0.410752 & 2.714842E-6 & 1.081073 & 1.081072 & 7.934216E-7
\\ 
$0.5$ & 0.521082 & 0.521095 & 1.326132E-5 & 1.127630 & 1.127625 & 4.774578E-6
\\ 
$0.6$ & 0.636604 & 0.636653 & 4.893639E-5 & 1.185485 & 1.185465 & 2.077300E-5
\\ 
$0.7$ & 0.758434 & 0.758583 & 1.490491E-4 & 1.255241 & 1.255169 & 7.230620E-5
\\ 
$0.8$ & 0.887710 & 0.888105 & 3.950285E-4 & 1.337648 & 1.337434 & 2.139083E-4
\\ 
$0.9$ & 0.025574 & 0.026516 & 9.426045E-4 & 1.433645 & 1.433086 & 5.592545E-4
\\ 
$1.0$ & 0.173128 & 0.175201 & 2.072716E-3 & 1.544407 & 1.543080 & 1.327116E-3
\\ \hline
\end{tabular}
}
\end{center}
\end{table}

\begin{figure}[tbp]
\centering
\includegraphics[width=3.50in]{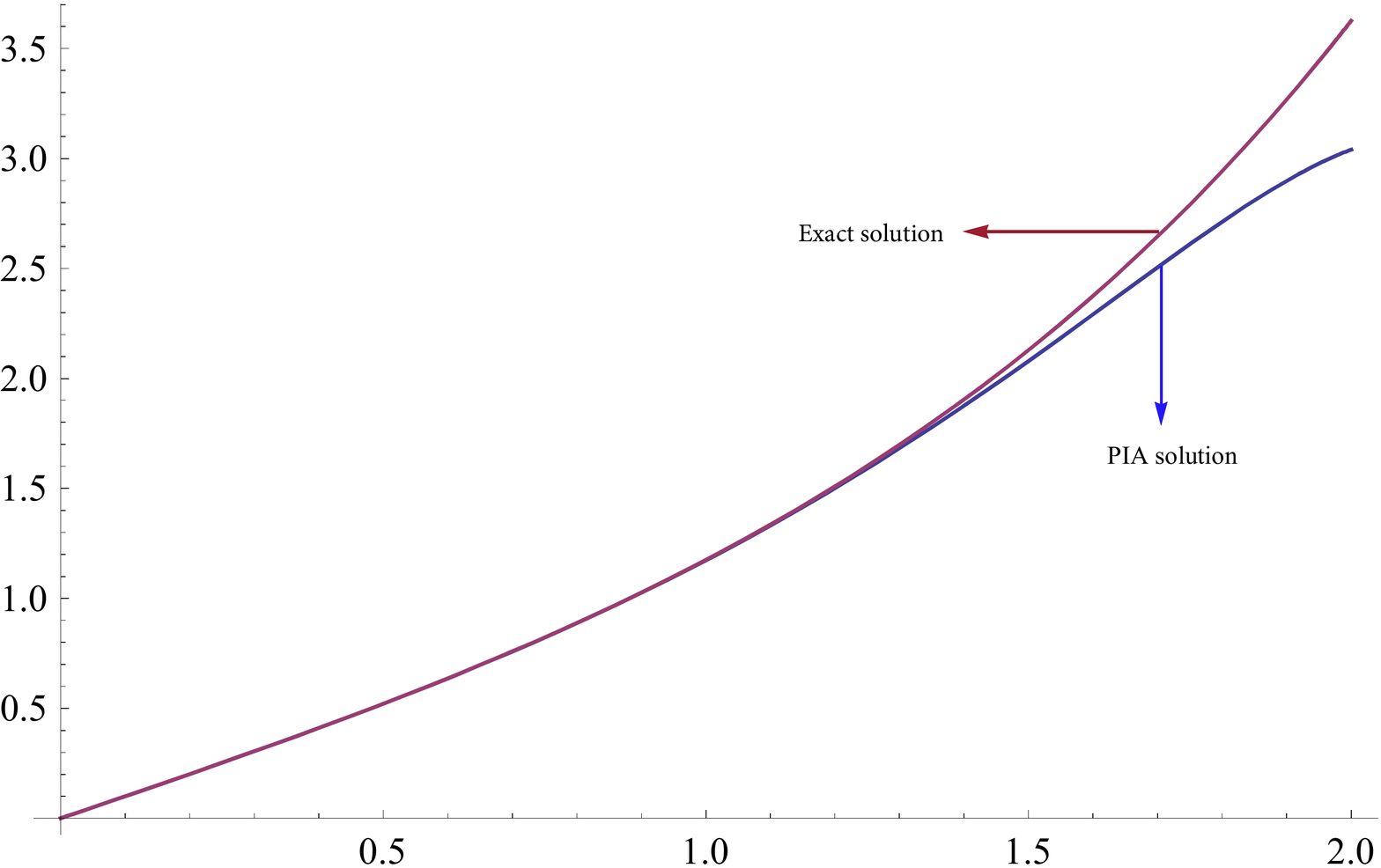}
\caption{Comparison of the PIA solution ($u_{3}(t)$) and exact solution for
Example 3.2. when $\alpha _{1} =\alpha _{2} =1$}
\end{figure}

\begin{figure}[tbp]
\centering
\includegraphics[width=3.50in]{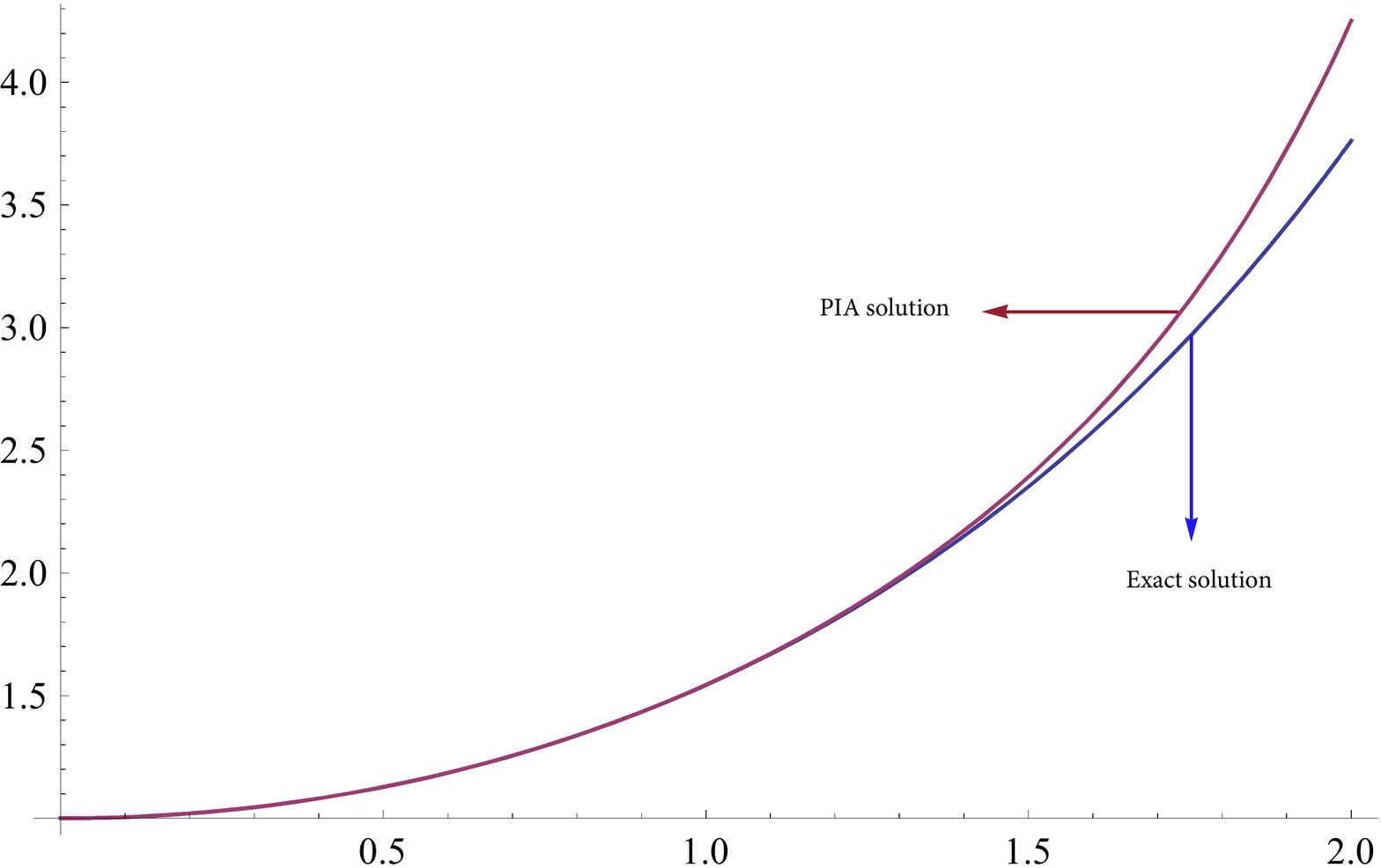}
\caption{Comparison of the PIA solution ($k_{3}(t)$) and exact solution for
Example 3.2. when $\alpha _{1} =\alpha _{2} =1$}
\end{figure}

\section{Conclusion}

\qquad In this study, Perturbation-Iteration Algorithm was introduced for
some Factional Differential Equations. It is clear that the method is very
simple and reliable perturbation-iteration technique and producing very approximate results. We expect that the present method could
used to calculate the approximate solutions of other types of fractional
differential equations.

\end{document}